\documentclass{amsart}
\usepackage{times,amssymb,array,mathrsfs}
\newcounter{theorem}

\newtheorem{Theorem}[theorem]{Theorem}
\newtheorem{Lemma}{Lemma}[section]
\newtheorem{Corollary}[theorem]{Corollary}
\newtheorem{Conjecture}[Lemma]{Conjecture}

\theoremstyle{definition}
\newtheorem{Example}[Lemma]{Example}
\newtheorem{Remark}[Lemma]{Remark}

\def\F{\mathbb F}

\def\Z{\mathbb Z}

\def\S{\mathcal S}

\def\SL{\mathrm{SL}}

\begin{document}

\title{Congruences of the partition function}
\author{Yifan Yang}
\address{Department of Applied Mathematics, National Chiao Tung University,
  1001 Ta Hsueh Road, Hsinchu, Taiwan 300}
\email{yfyang@math.nctu.edu.tw}
\date{\today}
\dedicatory{Dedicated to Professor B. C. Berndt on the occasion of his
  70th birthday}

\thanks{}

\begin{abstract} Let $p(n)$ denote the partition function. In this
  article, we will show that congruences of the form
  $$
    p(m^j\ell^kn+B)\equiv 0\mod m\ \text{ for all }\ n\ge 0
  $$
  exist for all primes $m$ and $\ell$ satisfying $m\ge 13$ and
  $\ell\neq 2,3,m$. Here the integer $k$ depends on the Hecke
  eigenvalues of a certain invariant subspace of
  $S_{m/2-1}(\Gamma_0(576),\chi_{12})$ and can be explicitly
  computed.

  More generally, we will show that for each integer $i>0$ there
  exists an integer $k$ such that for every non-negative
  integers $j\ge i$ with a properly chosen $B$ the congruence
  $$
    p(m^j\ell^kn+B)\equiv 0\mod m^i
  $$
  holds for all integers $n$ not divisible by $\ell$.
\end{abstract}

\maketitle

\begin{section}{Introduction} Let $p(n)$ denote the number of ways to
  write a positive integer $n$ as sums of positive integers. For
  convenience, we also set $p(0)=1$, $p(n)=0$ for $n<0$, and
  $p(\alpha)=0$ if $\alpha\not\in\Z$. A remarkable discovery of
  Ramanujan \cite{Ramanujan-Collected} is that the partition function
  $p(n)$ satisfies the congruences
  \begin{equation} \label{equation: Ramanujan congruences}
    p(An+B)\equiv 0 \mod m,
  \end{equation}
  for all non-negative integers $n$ for the triples
  $$
    (A,B,m)=(5,4,5), \ (7,5,7), \ (11,6,11).
  $$
  Ramanujan also conjectured that congruences \eqref{equation:
    Ramanujan congruences} exist for the cases $A=5^j$, $7^j$, or $11^j$.
%  $A$ and $m$ equal to higher powers of $5$, $7$, and $11$.
% hold for $(A,B,m)=(5^j,\delta_{5,j},5^j)$, $(7^j,\delta_{7,j},7^j)$,
%{\lfloor j/2\rfloor+1})$, and $(11^j,\delta_{11,j},11^j)$, where
%  $\delta_{m,j}$ denotes the unique integer satisfying
%  $24\delta_{m,j}\equiv 1\mod m^j$ and
%  $0<\delta_{m,j}<m^j$. This
  This conjecture was proved by Watson \cite{Watson-Crelle} for the
  cases of powers of $5$ and $7$ and Atkin \cite{Atkin-Glasgow}
  for the cases of powers of $11$. Since then, the problem of finding
  more examples of such congruences has attracted a great deal of
  attention. However, Ramanujan-type congruences appear to be very
  sparse. Prior to the late twentieth century, there are only a
  handful of such examples \cite{Atkin-PLMS,Atkin-OBrien-TAMS}. In those
  examples, the integers $A$ are no longer prime powers.

  It turns out that if we require the integer $A$ to be a prime,
  then the congruences proved or conjectured by Ramanujan are the only
  ones. This was proved recently in a remarkable paper of
  Ahlgren and Boylan \cite{Ahlgren-Boylan-Inventiones}. On the other
  hand, if $A$ is allowed to be a non-prime power, a surprising result
  of Ono \cite{Ono-Annals} shows that for each prime $m\ge 5$ and each
  positive integer $k$, a positive proportion of primes $\ell$ have
  the property
  \begin{equation} \label{equation: Ono congruences}
    p\left(\frac{m^k\ell^3n+1}{24}\right)\equiv 0\mod m
  \end{equation}
  for all non-negative integers $n$ relatively prime to $\ell$. This
  result was later extended to composite $m$, $(m,6)=1$, by Ahlgren
  \cite{Ahlgren-MAnn}. Neither of \cite{Ono-Annals} and
  \cite{Ahlgren-MAnn} addressed the algorithmic aspect of finding
  congruences of the form \eqref{equation: Ono congruences}. For the
  cases $m\in\{13,17,19,23,29,31\}$ this was done by
  Weaver \cite{Weaver-RamanujanJ}. In effect, she found 76,065 new
  congruences. For primes $m\ge 37$, this was addressed by Chua
  \cite{Chua-Archiv}. Although no explicit examples of congruences
  \eqref{equation: Ono congruences} for $m\ge 37$ were given in
  \cite{Chua-Archiv}, in principle, if one is patient enough, one will
  eventually find such congruences.
% Weaver's method was later refined by Chua
%  \cite{Chua-Archiv}. Potentially, Chua's
%  result can yield explicit examples of Ono's congruences
%  \eqref{equation: Ono congruences} for primes $m\ge 37$.

  Another remarkable discovery of Ono \cite[Theorem 5]{Ono-Annals} is
  that the partition function possesses certain periodic property
  modulo a prime $m$. Specifically, he showed that for every prime
  $m\ge 5$, there exist integers $0\le N(m)\le m^{48(m^3-2m+1)}$ and
  $1\le P(m)\le m^{48(m^3-2m+1)}$ such that
  \begin{equation} \label{equation: periodicity}
    p\left(\frac{m^in+1}{24}\right)\equiv
    p\left(\frac{m^{P(m)+i}n+1}{24}\right)\mod m
  \end{equation}
  for all non-negative integers $n$ and all $i\ge N(m)$. In
  \cite{Chua-Archiv}, Chua raised a conjecture (Conjecture
  \ref{conjecture: Chua} in Section \ref{section: main results}
  below), which, if is true, will greatly improve Ono's bound. (See
  Corollary \ref{corollary: periodicity} below.)
%  Then Chua's method
%  \cite{Chua-Archiv} improves the bound $m^{48(m^3-2m+1)}$ to
%  $m^{m-12\lfloor m/24\rfloor-1}$. We will briefly review the ideas of
%  Ono \cite{Ono-Annals}, Weaver \cite{Weaver-RamanujanJ}, and Chua
%  \cite{Chua-Archiv} in Section \ref{section: review}.

  In this note, we will obtain new congruences for the
  partition function and discuss related problems. In particular, we
  will show that there exist congruences of the form
  $$
    p(m^j\ell^kn+B)\equiv 0\mod m
  $$
  for \emph{all} primes $m$ and $\ell$ such that $m\ge 13$ and $\ell$
  not equal to $2,3,m$.

  \begin{Theorem} \label{theorem: simplified} Let
    $m$ and $\ell$ be primes such that $m\ge 13$ and $\ell\neq 2,3,m$.
    Then there exists an explicitly computable positive integer $k\ge
    2$ such that
    \begin{equation} \label{equation: Yang congruences}
      p\left(\frac{m^j\ell^{2k-1}n+1}{24}\right)\equiv 0\mod m
    \end{equation}
    for all non-negative integers $n$ relatively prime to $m$ and all
    positive integers $j$.
% The integer $k$ depends on the Hecke
%    operator $T_{\ell^2}$ on a certain Hecke invariant subspace of
%    $S_{m/2-1}(\Gamma_0(576),\chi_{12})$ and is explicitly computable.
%    $$
%      M=\begin{cases}m, &\text{if }m=13,17,19,23,29,31, \\
%       m^{\lfloor m/12\rfloor-\lfloor m/24\rfloor}+m-1,
%      &\text{if }m\ge 37. \end{cases}
%     (m^2-1)/2-1, &\text{if }m=37,41,43,47,53,59, \\
%     m^{3\lfloor m/12\rfloor-3\lfloor m/24\rfloor}/2, &\text{if }
%     m\ge 61.\end{cases}
%    $$
  \end{Theorem}

  For instance, in Section \ref{section: examples} we will find that
  for $m=37$ and arbitrary $j$, congruences \eqref{equation: Yang
    congruences} hold with
  $$ \extrarowheight3pt
  \begin{array}{c|ccccccccccccccc} \hline\hline
\ell & 5 & 7 & 11 & 13 & 17 & 19 & 23 & 29 & 31 & 41
     &43 &47 & 53 & 59 & 61 \\ \hline
k    &228&57 & 18 &684 & 38 & 38 &684 &684 &228 &171
     &18 &333& 18 & 12 & 684 \\ \hline\hline
  \end{array}
  $$
%  $$ \extrarowheight3pt
%  \begin{array}{c|ccccccccccccccccc} \hline\hline
%\ell & 5 & 7 & 11 & 17 & 19 & 23 & 29 & 31 & 37 & 41 & 43 & 47
%     & 53 & 59 & 61 & 67 & 73 \\ \hline
%k    &14 &14 & 14 & 7  & 14 & 3  &  6 & 12 & 14 & 12 &  7 & 12
%     & 7  &  2 & 13 & 12 & 2  \\
%\hline\hline
%  \end{array}
%  $$
  As far as we know, this is the first example in literature where a
  congruence \eqref{equation: Ramanujan congruences} modulo a prime
  $m\ge 37$ is explicitly given.

  Theorem \ref{theorem: simplified} is in fact a simplified version of
  one of the main results. (See Theorem \ref{theorem: congruences}).
  In the full version, we will see that the integer $k$ in Theorem
  \ref{theorem: simplified} can be determined quite explicitly in
  terms of the Hecke operators on a certain invariant subspace of the
  space $S_{m/2-1}(\Gamma_0(576),\chi_{12})$ of cusp forms of level
  $576$ and weight $m/2-1$ with character
  $\chi_{12}=\left(\frac{12}\cdot\right)$. To describe this invariant
  subspace and to see how it comes into play with congruences of the
  partition function, perhaps we should first
%
%  However, to see the relevance of our main results to the
%  study of congruences of the partition function, we will need to
  review the work of Ono \cite{Ono-Annals} and other subsequent
  papers \cite{Chua-Archiv,Weaver-RamanujanJ}. Thus, we will postpone
  giving the statements of our main results until Section
  \ref{section: main results}.

  Our method can be easily extended to obtain congruences of $p(n)$
  modulo a prime power. In Section \ref{section: generalizations},
  we will see that for each prime power $m^i$ and a prime $\ell\neq
  2,3,m$, there always exists a positive integer $k$ such that
  $$
    p\left(\frac{m^i\ell^{2k-1}n+1}{24}\right)\equiv 0\mod m^i
  $$
  for all positive integers $n$ not divisible by $\ell$. One example
  worked out in Section \ref{section: generalizations} is
  $$
    p\left(\frac{13^2\cdot5^{56783}n+1}{24}\right)\equiv 0\mod 13^2.
  $$
  In the same section, we will also discuss congruences of type
  $p(5^j\ell^kn+B)\equiv 0\mod 5^{j+1}$.
%  \eqref{equation: Ramanujan congruences} with $(A,m)=(5^j\ell^k,5^{j+1})$.
\smallskip

\noindent{\bf Notations.} Throughout the paper, we let
    $S_\lambda(\Gamma_0(N),\chi)$ denote the space of cusp forms of
    weight $\lambda$ and level $N$ with character $\chi$. By an
    invariant subspace of $S_\lambda(\Gamma_0(N),\chi)$ we mean a
    subspace that is invariant under the action of the Hecke algebra
    on the space.
% The notation $S_\lambda(\Gamma_0(N),\chi)_m$ means
%    the reduction of $\Z[[q]]\cap S_\lambda(\Gamma_0(N),\chi)$ modulo
%    a prime $m$.

    For a power series $f(q)=\sum a_f(n)q^n$ and a
    positive integer $N$, we let $U_N$ and $V_N$ denote the operators
    \begin{equation*}
    \begin{split}
      U_N&:f(q)\longmapsto f(q)\big|U_N:=\sum_{n=0}^\infty a_f(Nn)q^n, \\
      V_N&:f(q)\longmapsto f(q)\big|V_N:=\sum_{n=0}^\infty a_f(n)q^{Nn}.
    \end{split}
    \end{equation*}
    Moreover, if $\psi$ is a Dirichlet character, then $f\otimes\psi$
    denotes the twist $f\otimes\psi:=\sum a_f(n)\psi(n)q^n$.

    Finally, for a prime $m\ge 5$ and a positive integer $j$, we write
    $$
      F_{m,j}=\sum_{n\ge 0,m^jn\equiv-1\mod 24}
      p\left(\frac{m^jn+1}{24}\right)q^n.
    $$
\smallskip

  Note that we have
  \begin{equation} \label{equation: F|Um}
    F_{m,j}\big|U_m= F_{m,j+1}.
  \end{equation}
\end{section}

\begin{section}{Works of Ono \cite{Ono-Annals}, Weaver
    \cite{Weaver-RamanujanJ}, and Chua \cite{Chua-Archiv}}
  \label{section: review} In this section, we will review the ideas in
  \cite{Ono-Annals,Weaver-RamanujanJ,Chua-Archiv}.

  First of all, by a classical identity of Euler, we know that the
  generating function of $p(n)$ has an infinite product representation
  $$
    \sum_{n=0}^\infty p(n)q^n=\prod_{n=1}^\infty\frac1{1-q^n}.
  $$
  If we set $q=e^{2\pi i\tau}$, then we have
  $$
    q^{-1/24}\sum_{n=0}^\infty p(n)q^n=\eta(\tau)^{-1},
  $$
  where $\eta(\tau)$ is the Dedekind eta function. Now assume that $m$
  is a prime greater than $3$. Ono \cite{Ono-Annals} considered the
  function $\eta(m^k\tau)^{m^k}/\eta(\tau)$. On the one hand, one has
  $$
    \frac{\eta(m^k\tau)^{m^k}}{\eta(\tau)}\big|U_{m^k}
   =\prod_{n=1}^\infty(1-q^n)^{m^k}\cdot\left(\sum_{n=0}^\infty p(n)
    q^{n+(m^{2k}-1)/24}\right)\big| U_{m^k}.
  $$
  On the other hand, one has
  $$
    \frac{\eta(m^k\tau)^{m^k}}{\eta(\tau)}\equiv\eta(\tau)^{m^{2k}-1}
   =\Delta(\tau)^{(m^{2k}-1)/24} \mod m,
  $$
  where $\Delta(\tau)=\eta(\tau)^{24}$ is the normalized cusp from of
  weight $12$ on $\SL(2,\Z)$. From these, Ono
  \cite[Theorem 6]{Ono-Annals} deduced that
  $$
    F_{m,k}\equiv\frac{(\Delta(\tau)^{(m^{2k}-1)/24}\big|U_{m^k})\big|V_{24}}
    {\eta(24\tau)^{m^k}} \mod m.
  $$
  Now it can be verified that for $k=1$, the right-hand side of the
  above congruence is contained in the space
  $S_{(m^2-m-1)/2}(\Gamma_0(576m),\chi_{12})$ of cusp forms of level
  $576m$ and weight $(m^2-m-1)/2$ with character
  $\chi_{12}=\left(\frac{12}\cdot\right)$. Then by
  \eqref{equation: F|Um} and the fact that $U_m$ defines a linear
  map
  $$
    U_m:S_{\lambda+1/2}(\Gamma_0(4Nm),\psi)\to
        S_{\lambda+1/2}(\Gamma_0(4Nm),\psi\chi_m), \qquad
    \chi_m=\left(\frac m\cdot\right),
  $$
  one sees that
  $$
    F_{m,k}\equiv G_{m,k}=\sum a_{m,k}(n)q^n\mod m
  $$
  for some
%  Call the function on the right
%  $$
%    G_{m,k}:=\frac{(\Delta(\tau)^{(m^k-1)/24}\big|U_{m^k})\big|V_{24}}
%    {\eta(24\tau)^{m^k}}=\sum_{n=1}^\infty a_{m,k}(n)q^n.
%  $$
%  He then went on to show that $G_{m,k}$
%  $$
%    \frac{(\Delta(\tau)^{(m^k-1)/24}\big|U_{m^k})\big|V_{24}}
%    {\eta(24\tau)^{m^k}}
%  $$
  $G_{m,k}\in S_{(m^2-m-1)/2}(\Gamma_0(576m),\chi_{12}\chi_m^{k-1})$.
%, the space of
%  cusp forms of weight $(m^2-m-1)/2$ and level $576m$ with character
%  $\chi_{12}\chi^{k-1}$, where $\chi_{12}$ and $\chi$ are the Kronecker
%  characters associated with $\Q(\sqrt 3)$ and $\Q(\sqrt m)$,
%  respectively.

  Now recall the general Hecke theory for half-integral weight modular
  forms states that if $f(\tau)=\sum_{n=1}^\infty a_f(n)q^n\in
  S_{\lambda+1/2}(\Gamma_0(4N),\psi)$ and $\ell$ is a prime not
  dividing $4N$, then the Hecke operator defined by
  $$
    T_{\ell^2}:f(\tau)\mapsto\sum_{n=1}^\infty\left(
    a_f(\ell^2n)+\psi(\ell)\left(\frac{(-1)^\lambda n}{\ell}\right)
    \ell^{\lambda-1}a_f(n)+\psi(\ell^2)\ell^{2\lambda-1}a_f(n/\ell^2)\right)
    q^n
  $$
  sends $f(\tau)$ to a cusp form in the same space. In the situation
  under consideration, if $\ell$ is a prime not dividing $576m$ such
  that
  $$
    G_{m,k}\big|T_{\ell^2}\equiv 0\mod m,
  $$
  then we have
  \begin{equation*}
  \begin{split}
    0&\equiv(G_{m,k}\big|T_{\ell^2})\big|U_\ell \mod m \\
     &=\sum_{n=1}^\infty\left(a_{m,k}(\ell^3n)+\psi(\ell^2)
       \ell^{m^2-m-3}a_{m,k}(n/\ell)\right)q^n
  \end{split}
  \end{equation*}
  since $\left(\frac{\ell n}\ell\right)=0$. In particular, if $n$ is
  not divisible by $\ell$, then
  $$
    a_{m,k}(\ell^3n)\equiv 0\mod m,
  $$
  which implies
  $$
    p\left(\frac{m^k\ell^3n+1}{24}\right)\equiv 0\mod m.
  $$

  Finally, to show that there is a positive proportion of primes
  $\ell$ such that $G_{m,k}\big|T_{\ell^2}\equiv 0\mod m$, Ono invoked
  the Shimura correspondence between half-integral weight modular
  forms and integral weight modular forms \cite{Shimura-Annals} and a
  result of Serre \cite[6.4]{Serre}.

  As mentioned earlier, Ono \cite{Ono-Annals} did not address the
  issue of finding explicit congruences of the form \eqref{equation: Ono
    congruences}.
% For example, the only new congruence he gave is
%  $$
%    p\left(\frac{13\cdot59^3n+1}{24}\right)\equiv 0\mod 13,
%    \quad 59\nmid n.
%  $$
  However, Section 4 of \cite{Ono-Annals} did give us some hints on
  how one might proceed to discover new congruences, at least for
  small primes $m$. The key observation is the following.

  The modular form $G_{m,k}$ itself is in a vector space of big
  dimension, so to determine whether $G_{m,k}\big|T_{\ell^2}$ vanishes
  modulo $m$, one needs to compute the Fourier coefficients of
  $G_{m,k}$ for a huge number of terms.
% has a large weight, so potentially
%  a closed-form representation of $G_{m,k}$ (in terms of some familiar
%  modular forms such as Dedekind eta functions and Eisenstein series)
%  can be very complicated.
  However, it turns out that $F_{m,k}$ is
  congruent to another half-integral weight modular form of a much
  smaller weight. For example, using Sturm's theorem \cite{Sturm} Ono
  verified that
  \begin{equation} \label{equation: small weight}
  \begin{split}
    F_{13,2k+1}\equiv G_{13,2k+1}&\equiv 11\cdot 6^k\eta(24\tau)^{11} \mod 13, \\
    F_{13,2k+2}\equiv G_{13,2k+2}&\equiv 10\cdot 6^k\eta(24\tau)^{23} \mod 13
  \end{split}
  \end{equation}
  for all non-negative integers $k$. The modular form
  $\eta(24\tau)^{11}$ is in fact a Hecke eigenform. (The modular form
  $\eta(24\tau)^{23}$ is also a Hecke eigenform as we shall see in
  Section \ref{section: main results}.) More generally, for
  $m\in\{13,17,19,23,29,31\}$, it is shown in \cite[Section
  4]{Ono-Annals}, \cite[Proposition 6]{Guo-Ono-IMRN} and
  \cite[Proposition 5]{Weaver-RamanujanJ} that $G_{m,1}$ is congruent
  to a Hecke eigenform of weight $m/2-1$. Using this observation,
  Weaver \cite{Weaver-RamanujanJ} then devised an algorithm to find
  explicit congruences of the form \eqref{equation: Ono congruences}
  for $m\in\{13,17,19,23,29,31\}$.

  The proof of congruences \eqref{equation: small weight}
  given in \cite{Guo-Ono-IMRN} and \cite{Weaver-RamanujanJ} is
  essentially ``verification'' in the sense that they all used Sturm's
  criterion \cite{Sturm}. That is, by Sturm's theorem to show that two
  modular forms on a congruence subgroup $\Gamma$ are congruent to
  each other modulo a prime $m$, it suffices to compare sufficiently
  many coefficients, depending on the weight and index
  $(\SL(2,\Z):\Gamma)$. Naturally, this kind of argument will not be
  very useful in proving general results. In
  \cite{Chua-Archiv}, Chua found a more direct way to prove
  congruences \eqref{equation: small weight} for $F_{m,1}$. In
  particular, he \cite[Theorem 1.1]{Chua-Archiv} was able to show that
  for each prime $m\ge 5$, $F_{m,1}$ is congruent to a modular form of
  weight $m/2-1$ modulo $m$.

  Instead of the congruence
  $$
    \frac{\eta(m\tau)^m}{\eta(\tau)}\equiv\eta(\tau)^{m^2-1} \mod m
  $$
  used by Ono, Chua considered the congruence
  $$
    \frac{\eta(m\tau)^m}{\eta(\tau)}\equiv\eta(m\tau)^{m-1}\eta(\tau)^{m-1}
    \mod m
  $$
  as the starting point. The function on the right is a modular form
  of weight $m-1$ on $\Gamma_0(m)$. Thus, by the level reduction lemma
  of Atkin and Lehner \cite[Lemma 7]{Atkin-Lehner-MAnn}, one has
  $$
    \eta(m\tau)^{m-1}\eta(\tau)^{m-1}\big|(U_m+m^{(m-1)/2-1}W_m) \in
    S_{m-1}(\SL(2,\Z)),
  $$
  where $W_m$ denotes the Atkin-Lehner involution. It follows that
  $$
    F_{m,1}=\frac1{\eta(24\tau)}\big|U_m
    \equiv\frac{f_m(24\tau)}{\eta(24\tau)^m} \mod m
  $$
  for some cusp form $f_m(\tau)\in S_{m-1}(\SL(2,\Z))$.
  (Incidently, this also proves Ramanujan's congruences for
  $m=5,7,11$, since there are no non-trivial cusp forms of weight
  $4,6,10$.) By examining the order of vanishing of $f_m(\tau)$ at
  $\infty$, Chua \cite[Theorem 1.1]{Chua-Archiv} then concluded that
  if we let $r_m$ denote the integer in the range $0<r_m<24$ such that
  $m\equiv-r_m\mod 24$, then
  $$
    F_{m,1}\equiv\eta(24\tau)^{r_m}\phi_m(24\tau)
  $$
  for some modular form $\phi_m$ on $\SL(2,\Z)$ of weight
  $(m-r_m-2)/2$. Furthermore, based on an extensive numerical
  computation, Chua made the following conjecture.

\begin{Conjecture}[{Chua \cite[Conjecture 1]{Chua-Archiv}}]
  \label{conjecture: Chua} Let $m\ge 13$ be a prime and $r_m$ be the
  integer in the range $0<r_m<24$ such that $m\equiv-r_m\mod 24$. Set
  $$
    r_{m,j}=\begin{cases}r_m, &\text{if }j\text{ is odd}, \\
    23, &\text{if }j\text{ is even}. \end{cases}
  $$
  Then
  $$
    F_{m,j}\equiv\eta(24\tau)^{r_{m,j}}\phi_{m,j}(24\tau) \mod m
  $$
  for some modular form $\phi_{m,j}(\tau)$ on $\SL(2,\Z)$, where the
  weight of $\phi_{m,j}$ is $(m-r_m-2)/2$ if $j$ is odd and is $m-13$
  if $j$ is even.
\end{Conjecture}

In \cite[Section 4]{Chua-Archiv}, Chua established the induction step
for the case of even $j$ assuming the conjecture holds for odd $j-1$.
However, as remarked by Chua, it appears difficult to prove the
induction step from cases of even $j-1$ to cases of odd $j$. In the
next section, we will see that this conjecture is a simple consequence
of our Theorem \ref{theorem: invariant subspace}.

\begin{Remark} Professor H. H. Chan has kindly informed us that Serre
  has indicated to him an argument to establish Conjecture
  \ref{conjecture: Chua}. The argument will be given in a
  forthcoming article \cite{Chan-Serre}.
\end{Remark}
\end{section}

\begin{section}{Statements of main results} \label{section: main results}

  The functions $\eta(24\tau)^{r_{m,k}}\phi_{m,k}(24\tau)$ appearing
  in Chua's conjecture (Conjecture \ref{conjecture: Chua}) are all
  half-integral weight modular forms of level $576$ and character
  $\chi_{12}$. Thus, our first main result is concerned with the space
  $S_{\lambda+1/2}(\Gamma_0(576),\chi_{12})$.

  \begin{Theorem} \label{theorem: invariant subspace} Let $r$ be an
    odd integer with $0<r<24$. Let $s$ be a non-negative even integer.
    Then the space
    \begin{equation} \label{equation: Srs}
      \S_{r,s}:=\left\{\eta(24\tau)^rf(24\tau):f(\tau)
      \in M_s(\SL(2,\Z))\right\}
    \end{equation}
    is an invariant subspace of $S_{s+r/2}(\Gamma_0(576),\chi_{12})$
    under the action of the Hecke algebra. That is, for all primes
    $\ell\neq 2,3$ and all $f\in\S_{r,s}$, we have
    $f\big|T_{\ell^2}\in\S_{r,s}$.
  \end{Theorem}

  The following corollary is immediate.

  \begin{Corollary} Let $r$ be an odd integer with $0<r<24$. Let
    $E_4(\tau)$ and $E_6(\tau)$ be the Eisenstein series of weights
    $4$ and $6$ on $\SL(2,\Z)$ and $f(\tau)$ be one of the function
    $1$, $E_4(\tau)$, $E_6(\tau)$, $E_4(\tau)^2$,
    $E_4(\tau)E_6(\tau)$, and $E_4(\tau)^2E_6(\tau)$. Then the
    function $\eta(24\tau)^rf(24\tau)$ is a Hecke eigenform.
    In particular, for $m\in\{13,17,19,23,29,31\}$, the function
    $$
      \eta(24\tau)^{r_m}\phi_{m,1}(24\tau)
    $$
    in Conjecture \ref{conjecture: Chua} is a Hecke eigenform.
  \end{Corollary}

  Note that the assertion about $g_m:=\eta(24\tau)^{r_m}\phi_{m,1}$
  was already proved in Proposition 6 of \cite{Guo-Ono-IMRN}. In the
  same proposition, it was also proved that the image of $g_m$ under
  the Shimura correspondence is $G_m\otimes\chi_{12}$, where $G_m$ is
  the unique normalized newform of weight $m-3$ on $\Gamma_0(6)$ whose
  eigenvalues for the Atkin-Lehner involutions $W_2$ and $W_3$ are
  $-\left(\frac2m\right)$ and $-\left(\frac3m\right)$, respectively.

  We now apply Theorem \ref{theorem: invariant subspace} to study
  congruences of the partition function. We first consider Conjecture
  \ref{conjecture: Chua}. Observe that the Hecke
  operator $T_{m^2}$ is the same as the operator $U_{m^2}$ modulo $m$.
  Also, the case $j=1$ and the induction step from $j=1$ to $j=2$
  have already been proved in \cite{Chua-Archiv}. Thus, from Theorem
  \ref{theorem: invariant subspace} we immediately conclude that
  Chua's conjecture indeed holds in general.

  \begin{Corollary}[Conjecture of Chua] \label{corollary: Chua}
  Let $m\ge 13$ be a prime and $r_m$ be the integer in the range
  $0<r_m<24$ such that $m\equiv-r_m\mod 24$. Set
  $$
    r_{m,j}=\begin{cases}r_m, &\text{if }j\text{ is odd}, \\
    23, &\text{if }j\text{ is even}. \end{cases}
  $$
  Then
  $$
    F_{m,j}\equiv\eta(24\tau)^{r_{m,j}}\phi_{m,j}(24\tau) \mod m
  $$
  for some modular form $\phi_{m,j}(\tau)$ on $\SL(2,\Z)$, where the
  weight of $\phi_{m,j}$ is $(m-r_m-2)/2$ if $j$ is odd and is $m-13$
  if $j$ is even.
  \end{Corollary}

  \begin{Remark} Note that for odd $j$, we have
  \begin{equation} \label{equation: dimension}
    \dim\S_{r_m,(m-r_m-2)/2}=\left\lfloor\frac m{12}\right\rfloor
   -\left\lfloor\frac m{24}\right\rfloor.
  \end{equation}
  To see this, we observe that $\dim
  M_\lambda(\SL(2,\Z))-\lfloor\lambda/12\rfloor$ is periodic of
  period $12$. Thus, to show \eqref{equation: dimension}, we only
  need to verify case by case according the residue of $m$ modulo $24$.
  \end{Remark}

  Using the pigeonhole principle, one can see that Theorem
  \ref{theorem: invariant subspace} also yields Ono's periodicity
  result \eqref{equation: periodicity}, with an improved bound.

  \begin{Corollary} \label{corollary: periodicity} Let $m\ge 5$ be a
    prime. Then there exist integers $0\le N(m)\le m^{A(m)}$ and $0\le
    P(m)\le m^{A(m)}$ such that
  $$
    p\left(\frac{m^in+1}{24}\right)\equiv
    p\left(\frac{m^{P(m)+i}n+1}{24}\right)\mod m
  $$
  for all non-negative integers $n$, where
  \begin{equation} \label{equation: A(m)}
    A(m)=2\dim M_{(m-r_m-2)/2}(\SL(2,\Z))
  \end{equation}
  and $r_m$ is the integer satisfying $0<r_m<24$ and $m\equiv-r_m\mod 24$.
  \end{Corollary}

%  In fact, without much difficulty, we can also obtain another similar
%  periodicity result from Theorem \ref{theorem: invariant subspace}.

%  \begin{Corollary} Let $m\ge 5$ be a prime and $\ell\ge 5$ be another
%    prime different from $m$. Then there exist integers $0\le N(m,\ell)\le
%    m^{2A(m)}+1$ and $0\le P(m,\ell)\le m^{2A(m)}$ such that
%  $$
%    p\left(\frac{m\ell ^in+1}{24}\right)\equiv
%    p\left(\frac{m\ell^{P(m)+i}n+1}{24}\right)\mod m
%  $$
%  for all non-negative integers $n$ and all integers $i\ge N(m,\ell)$,
%  where $A(m)$ is the same as \eqref{equation: A(m)}.
%  \end{Corollary}

  \begin{Corollary} \label{corollary: recursion} Let $r$ be an odd
    integer satisfying $0<r<24$ and $s$ be a non-negative even
    integer. Let $\S_{r,s}$ be defined as \eqref{equation: Srs} and
    $\{f_1,\ldots,f_t\}$ be a $\Z$-basis for the $\Z$-module
    $\Z[[q]]\cap\S_{r,s}$. Given a prime $\ell\ge 5$, assume that
    $A$ is the $t\times t$ matrix such that
    $$
      \begin{pmatrix}f_1\\\vdots\\f_t\end{pmatrix}\big|T_{\ell^2}
     =A\begin{pmatrix}f_1\\\vdots\\f_t\end{pmatrix}.
    $$
    Then we have
    $$
      \begin{pmatrix}f_1\\\vdots\\f_t\end{pmatrix}\big|U_{\ell^2}^k
  =A_k\begin{pmatrix}f_1\\\vdots\\f_t\end{pmatrix}
  +B_k\begin{pmatrix}g_1\\\vdots\\g_t\end{pmatrix}
  +C_k\begin{pmatrix}f_1\\\vdots\\f_t\end{pmatrix}\big|V_{\ell^2},
    $$
    where $g_j=f_j\otimes \left(\frac\cdot\ell\right)$, and $A_k$,
    $B_k$, and $C_k$ are $t\times t$ matrices satisfying
  $$
    \begin{pmatrix}A_k\\A_{k-1}\end{pmatrix}
   =\begin{pmatrix}A&-\ell^{r+2s-2}I_t\\I_t&0\end{pmatrix}^k
    \begin{pmatrix}I_t\\0\end{pmatrix},
  $$
  and
  $$
    B_k=-\ell^{s+(r-3)/2}\left(\frac{(-1)^{(r-1)/2}12}\ell\right)A_{k-1},
    \qquad C_k=-\ell^{r+2s-2}A_{k-1}.
  $$
%    the recursive relations
%    $$
%      A_1=A, \quad
%      B_1=-\ell^{s+(r-3)/2}\left(\frac{(-1)^{(r-1)/2}12}\ell\right)I_t,
%      \quad
%      C_1=-\ell^{r+2s-2}I_t,
%    $$
%    and
%    $$
%      A_k=A_{k-1}A_1+C_{k-1}, \qquad B_k=A_{k-1}B_1, \qquad
%      C_k=A_{k-1}C_1.
%    $$
  \end{Corollary}

%  \begin{Remark} Note that $A_k$, $B_k$, and $C_k$ are all polynomials
%    in $A$. Thus, they commute with each other.
%  \end{Remark}

  \begin{Theorem} \label{theorem: congruences} Let $m\ge 13$ be a
    prime and $j$ be a positive integer. Set $r_m$ to be the integer
    satisfying $0<r_m<24$ and $m\equiv-r_m\mod 24$. Let
    $$
      t=\left\lfloor\frac m{12}\right\rfloor
       -\left\lfloor\frac m{24}\right\rfloor
    $$
    be the dimension of $\S_{r_m,(m-r_m-2)/2}$ and assume that
    $\{f_1,\ldots,f_t\}$ is a $\Z$-basis for the $\Z$-module
    $\Z[[q]]\cap\S_{r_m,(m-r_m-2)/2}$.
    Let $\ell$ be a prime different from $2$, $3$, and $m$, and assume
    that $A$ is the $t\times t$ matrix such that
    $$
      \begin{pmatrix}f_1\\\vdots\\f_t\end{pmatrix}\big|T_{\ell^2}
     =A\begin{pmatrix}f_1\\\vdots\\f_t\end{pmatrix}.
    $$
    Assume that the order of the square matrix
    \begin{equation} \label{equation: matrix}
      \begin{pmatrix} A&-\ell^{m-4}I_t\\I_t&0\end{pmatrix} \mod m
    \end{equation}
    in $\mathrm{PGL}(2t,\mathbb F_m)$ is $K$, then we have
    \begin{equation} \label{equation: result}
      p\left(\frac{m^j\ell^{2uK-1}+1}{24}\right)\equiv 0\mod m
    \end{equation}
    for all positive integers $j$ and $u$ and all positive integers $n$ not
    divisible by $\ell$.

    Also, if the order of the matrix \eqref{equation: matrix} in
    $\mathrm{GL}(2t,\F_m)$ is $M$, then we have
    \begin{equation} \label{equation: periodicity result}
      p\left(\frac{m^j\ell^i n+1}{24}\right)\equiv
      p\left(\frac{m^j\ell^{2M+i}n+1}{24}\right) \mod m
    \end{equation}
    for all non-negative integer $i$ and all positive integers $j$ and
    $n$.
  \end{Theorem}

  \begin{Remark} Note that if the matrix $A$ in the above theorem
    vanishing modulo $m$, then the matrix in \eqref{equation: matrix}
    has order $2$ in $\mathrm{PGL}(2t,\F_m)$, and the conclusion of
    the theorem asserts that
    $$
      p\left(\frac{m^j\ell^3n+1}{24}\right)\equiv 0\mod m.
    $$
    This is the congruence appearing in Ono's theorem.
  \end{Remark}

  \begin{Remark} In general, the integer $K$ in Theorem \ref{theorem:
      congruences} may not be the smallest positive integer such that
    congruence \eqref{equation: Yang congruences} holds. We choose to
    state the theorem in the current form because of its simplicity.
    See the remark following the proof of Theorem \ref{theorem:
      congruences}.
  \end{Remark}
\end{section}

\begin{section}{Proof of Theorem \ref{theorem: invariant subspace}}
  We first recall the following lemma of Atkin and Lehner.

  \begin{Lemma}[{\cite[Lemma 7]{Atkin-Lehner-MAnn}}]
    \label{lemma: Atkin-Lehner} 
  Let $f$ be a cusp form of weight $s$ on $\Gamma_0(N)$ and $\ell$ be
  a prime. Then
  \begin{enumerate}
  \item[(a)] If $\ell|N$, then $f\big|U_\ell$ is a cusp form on
    $\Gamma_0(N)$. Furthermore, if $\ell^2|N$, then $f\big|U_\ell$ is
    modular on $\Gamma_0(N/\ell)$.
  \item[(b)] If $\ell|N$ but $\ell^2\nmid N$, then
    $f\big|(U_\ell+\ell^{s/2-1}W_\ell)$ is a cusp form on
    $\Gamma_0(N/\ell)$.
  \end{enumerate}
  \end{Lemma}

  The following transformation formula for $\eta(\tau)$ is frequently
  used.

  \begin{Lemma}[{\cite[pp.125--127]{Weber}}] \label{lemma: eta} For
  $$
    \gamma=\begin{pmatrix}a&b\\ c&d\end{pmatrix}\in\SL_2(\mathbb Z),
  $$
  the transformation formula for $\eta(\tau)$ is given by, for $c=0$,
  $$
    \eta(\tau+b)=e^{\pi ib/12}\eta(\tau),
  $$
  and, for $c>0$,
  $$
    \eta(\gamma\tau)=\epsilon(a,b,c,d)\sqrt{\frac{c\tau+d}i}\eta(\tau)
  $$
  with
  \begin{equation*}
%  \label{epsilon1}
    \epsilon(a,b,c,d)=
    \begin{cases}\displaystyle
      \left(\frac dc\right)i^{(1-c)/2}
      e^{\pi i\left(bd(1-c^2)+c(a+d)\right)/12},
      &\text{if }c\text{ is odd},\\
     \displaystyle 
     \left(\frac cd\right)e^{\pi i\left(ac(1-d^2)+d(b-c+3)\right)/12},
      &\text{if }d\text{ is odd},
    \end{cases}
  \end{equation*}
  where $\displaystyle\left(\frac dc\right)$ is the Legendre-Jacobi symbol.
  \end{Lemma}

  We now prove Theorem \ref{theorem: invariant subspace}.
  Assume that $g(\tau)\in\S_{r,s}$, say
  $g(\tau)=\eta(24\tau)^rf(24\tau)$ for some $f(\tau)\in
  M_s(\SL(2,\Z))$. Assume $g(\tau)=q^r\sum_{n=0}^\infty a(n)q^{24n}$.
  Then by the definition of $T_{\ell^2}$, we have
  \begin{equation} \label{equation: theorem 2 1}
  \begin{split}
    g\big|T_{\ell^2}
    &=\sum_{n\ge 0,n\equiv-r/24\mod\ell^2}a(n)q^{(24n+r)/\ell^2} \\
    &\quad\qquad+\ell^{s+(r-3)/2}
     \left(\frac{(-1)^{(r-1)/2}12}\ell\right)
     \sum_{n=0}^\infty\left(\frac{24n+r}\ell\right)a(n)q^{24n+r} \\
    &\quad\qquad+\ell^{r+2s-2}\sum_{n=0}^\infty a(n)q^{\ell^2(24n+r)}.
  \end{split}
  \end{equation}
  We now consider the function
  $$
    h(\tau)=\eta(\ell^2\tau)^{24-r}g(\tau/24)
    =\eta(\ell^2\tau)^{24-r}\eta(\tau)^rf(\tau).
  $$
  Using Newman's criterion \cite[Theorem 1]{Newman-PLMS}, we see that
  $h(\tau)$ is a cusp form of weight $s+12$ on $\Gamma_0(\ell^2)$.
  Then by Lemma \ref{lemma: Atkin-Lehner},
  $h\big|(U_{\ell^2}+\ell^{s/2+5}U_\ell W_\ell)$
  is a cusp form on $\SL(2,\Z)$, that is,
  \begin{equation} \label{equation: theorem 2 2}
    h\big|(U_{\ell^2}+\ell^{s/2+5}U_\ell W_\ell)=\eta(\tau)^{24}
    \tilde h(\tau)
  \end{equation}
  for some modular form $\tilde h(\tau)$ of weight $s$ on $\SL(2,\Z)$.
  We claim that
  \begin{equation} \label{equation: theorem 2 3}
    h\big|(U_{\ell^2}+\ell^{s/2+5}U_\ell W_\ell)=\eta(\tau)^{24-r}
    (g\big|T_{\ell^2})(\tau/24).
  \end{equation}
  Once this is proved, by comparing \eqref{equation: theorem 2 2} and
  \eqref{equation: theorem 2 3}, we immediately get Theorem
  \ref{theorem: invariant subspace}. We now verify \eqref{equation:
    theorem 2 3}.

  By the definition of $U_{\ell^2}$ we have
  \begin{equation} \label{equation: theorem 2 4}
  \begin{split}
    h\big|U_{\ell^2}&=\left(\prod_{n=1}^\infty(1-q^{\ell^2n})^{24-r}
    \sum_{n=0}^\infty
    a(n)q^{\ell^2-r(\ell^2-1)/24+n}\right)\big|U_{\ell^2} \\
    &=q\prod_{n=1}^\infty(1-q^n)^{24-r}
      \sum_{n\ge
        0,n\equiv-r/24\mod\ell^2}a(n)q^{(24n-r(\ell^2-1))/24\ell^2} \\
    &=\eta(\tau)^{24-r}\sum_{n\ge 0,n\equiv-r/24\mod\ell^2}
      a(n)q^{(24n+r)/24\ell^2}.
  \end{split}
  \end{equation}
  The term involving $U_\ell W_\ell$ is more complicated. We have
  \begin{equation*}
  \begin{split}
    h\big|U_\ell W_\ell
   =\left(\frac1\ell\sum_{k=0}^{\ell-1}h\big|
    \begin{pmatrix}1&k\\0&\ell\end{pmatrix}\right)\big|W_\ell
  =\ell^{-s/2-7}\tau^{-s-12}\sum_{k=0}^{\ell-1}h\big|
    \begin{pmatrix}1&k\\0&\ell\end{pmatrix}
    \begin{pmatrix}0&-1\\\ell&0\end{pmatrix}.
  \end{split}
  \end{equation*}
  The term $k=0$ gives us
  \begin{equation*}
    \ell^{-s/2-7}\tau^{-s-12}h(-1/\ell^2\tau)
   =\ell^{-s/2-7}\tau^{-s-12}\eta(-1/\tau)^{24-r}\eta(-1/\ell^2\tau)^r
    f(-1/\ell^2\tau).
  \end{equation*}
  Using the formula $\eta(-1/\tau)=\sqrt{(\tau/i)}\eta(\tau)$ and the
  assumption that $f(\tau)\in M_s(\SL(2,\Z))$, this reduces to
  \begin{equation} \label{equation: k=0}
    \ell^{3s/2+r-7}\eta(\tau)^{24-r}\eta(\ell^2\tau)^rf(\ell^2\tau)
   =\ell^{3s/2+r-7}\eta(\tau)^{24-r}\sum_{n=0}^\infty a(n)q^{\ell^2(24n+r)/24}.
  \end{equation}
  We now consider the contribution from the cases $k\neq 0$.

  We have
  $$
    \eta(\ell^2\tau)\big|\begin{pmatrix}k\ell&-1\\\ell^2&0
    \end{pmatrix}
   =\eta\left(\frac{k\ell\tau-1}\tau\right).
  $$
  By Lemma \ref{lemma: eta}, this is equal to
  \begin{equation} \label{equation: k>0 eta(ell2)}
    \eta(\ell^2\tau)\big|\begin{pmatrix}k\ell&-1\\\ell^2&0
    \end{pmatrix}
   =e^{2\pi ik\ell/24}\sqrt{\frac\tau i}\eta(\tau).
  \end{equation}
  For $\eta(\tau)$ and $f(\tau)$, we observe that
  $$
    \frac{k\ell\tau-1}{\ell^2\tau}
   =\begin{pmatrix}k&u\\\ell&k'\end{pmatrix}(\tau-k'/\ell),
  $$
  where $k'$ denotes the multiplicative inverse of $k$ modulo $\ell$
  and $u=(kk'-1)/\ell$. Thus, by Lemma \ref{lemma: eta},
  \begin{equation} \label{equation: k>0 eta}
    \eta(\tau)\big|\begin{pmatrix}k\ell&-1\\\ell^2&0\end{pmatrix}
   =\left(\frac{k'}\ell\right)i^{(1-\ell)/2}
    e^{2\pi i\ell(k+k')/24}\sqrt{\frac{\ell\tau}i}\eta
    \left(\tau-\frac{k'}\ell\right).
  \end{equation}
  Also,
  \begin{equation} \label{equation: k>0 f}
    f(\tau)\big|\begin{pmatrix}k\ell&-1\\\ell^2&0\end{pmatrix}
   =(\ell\tau)^sf\left(\tau-\frac{k'}\ell\right).
  \end{equation}
  Combining \eqref{equation: k>0 eta(ell2)}, \eqref{equation: k>0
    eta}, and \eqref{equation: k>0 f}, we obtain
  \begin{equation*}
  \begin{split}
   &\ell^{-s/2-7}\tau^{-s-12}h\big|
    \begin{pmatrix}k\ell&-1\\\ell^2&0\end{pmatrix} \\
   &\qquad=\ell^{(s+r)/2-7}\left(\frac{k'}\ell\right)
    i^{r(1-\ell)/2}e^{2\pi ir\ell k'/24}\eta(\tau)^{24-r}
    \eta\left(\tau-\frac{k'}\ell\right)^r
    f\left(\tau-\frac{k'}\ell\right),
  \end{split}
  \end{equation*}
  and
  \begin{equation} \label{equation: k>0 h 1}
  \begin{split}
   &\ell^{-s/2-7}\tau^{-s-12}\sum_{k=1}^{\ell-1}h\big|
    \begin{pmatrix}k\ell&-1\\\ell^2&0\end{pmatrix} \\
   &\qquad=\ell^{(s+r)/2-7}i^{r(1-\ell)/2}\eta(\tau)^{24-r}
    \sum_{k=1}^{\ell-1}\left(\frac{k}\ell\right)
    e^{2\pi ir\ell k/24}g\left(\frac{\tau-k/\ell}{24}\right).
  \end{split}
  \end{equation}
  The sum in the last expression is equal to
  \begin{equation} \label{equation: k>0 h 2}
  \begin{split}
    \sum_{k=1}^{\ell-1}\left(\frac{k}\ell\right)
    e^{2\pi irk(\ell^2-1)/24\ell}
    \sum_{n=0}^\infty e^{-2\pi ikn/\ell}a(n)q^{n+r/24}.
  \end{split}
  \end{equation}
  With the well-known evaluation
  $$
    \sum_{k=1}^{\ell-1}\left(\frac k\ell\right)
    e^{2\pi ikn/\ell}
   =\left(\frac n\ell\right)i^{(\ell-1)^2/4}\sqrt\ell
  $$
  of the Gaussian sum, \eqref{equation: k>0 h 2} can be simplified to
  \begin{equation*}
  \begin{split}
   &i^{(\ell-1)^2/4}\sqrt\ell\sum_{n=0}^\infty
    \left(\frac{r(\ell^2-1)/24-n}\ell\right)a(n)q^{n+r/24} \\
   &\qquad=i^{(\ell-1)^2/4}\sqrt\ell\left(\frac{-24}\ell\right)
    \sum_{n=0}^\infty \left(\frac{24n+r}\ell\right)a(n)
    q^{n+r/24}.
  \end{split}
  \end{equation*}
  Substituting this into \eqref{equation: k>0 h 1} and using
  $$
    \left(\frac{-1}\ell\right)=(-1)^{(\ell-1)/2}, \quad
    \left(\frac  2\ell\right)=(-1)^{(\ell^2-1)/8},
  $$
   we arrive at
  \begin{equation} \label{equation: k>0 h 3}
  \begin{split}
   &\ell^{-s/2-7}\tau^{-s-12}\sum_{k=1}^{\ell-1}h\big|
    \begin{pmatrix}k\ell&-1\\\ell^2&0\end{pmatrix} \\
   &\quad=\ell^{(s+r+1)/2-7}i^{r(1-\ell)/2+(\ell-1)^2/4}
    \left(\frac{-24}\ell\right)\eta(\tau)^{24-r}
    \sum_{n=0}^\infty\left(\frac{24n+r}\ell\right)a(n)q^{n+r/24} \\
   &\quad=\ell^{(s+r+1)/2-7}\left(\frac{-1}\ell\right)^{(r-1)/2}
    \left(\frac{12}\ell\right)\eta(\tau)^{24-r}
    \sum_{n=0}^\infty\left(\frac{24n+r}\ell\right)a(n)q^{n+r/24}.
  \end{split}
  \end{equation}
  Together with \eqref{equation: k=0}, \eqref{equation: k>0 h 3}
  implies that
  \begin{equation} \label{equation: k>0 h last}
  \begin{split}
   &\ell^{s/2+5}h\big|U_\ell W_\ell
   =\ell^{2s+r-2}\eta(\tau)^{24-r}\sum_{n=0}^\infty
    a(n)q^{\ell^2(24n+r)/24} \\
   &\quad\qquad+\ell^{s+(r-3)/2}\eta(\tau)^{24-r}
    \left(\frac{(-1)^{(r-1)/2}12}\ell\right)
    \sum_{n=0}^\infty\left(\frac{24n+r}\ell\right)a(n)q^{n+r/24}.
  \end{split}
  \end{equation}
  Comparing \eqref{equation: theorem 2 4} and \eqref{equation: k>0 h
    last} with \eqref{equation: theorem 2 1}, we see that
  \eqref{equation: theorem 2 3} indeed holds. The proof of Theorem
  \ref{theorem: invariant subspace} is now complete.
\end{section}

\begin{section}{Proof of Corollary \ref{corollary: recursion} and
  Theorem \ref{theorem: congruences}} \label{section: proof 2}
\begin{proof}[Proof of Corollary \ref{corollary: recursion}] By the
  definition of $T_{\ell^2}$, we have
  $$
    \begin{pmatrix}f_1\\\vdots\\ f_t\end{pmatrix}\big|U_{\ell^2}
   =A_1\begin{pmatrix}f_1\\\vdots\\ f_t\end{pmatrix}
   +B_1\begin{pmatrix}g_1\\\vdots\\g_t\end{pmatrix}
   +C_1\begin{pmatrix}f_1\\\vdots\\
    f_t\end{pmatrix}\big|V_{\ell^2},
  $$
  where $g_t=f_t\otimes\left(\frac\cdot\ell\right)$ and
  $$
    A_1=A, \quad
    B_1=-\ell^{s+(r-3)/2}\left(\frac{(-1)^{(r-1)/2}12}\ell\right)I_t,
    \quad
    C_1=-\ell^{r+2s-2}I_t.
  $$
  Now we make the key observation
  $$
    g_j\big|U_{\ell^2}=0, \qquad f_j\big|V_{\ell^2}\big|U_{\ell^2}=f_j,
  $$
  from which we obtain
  $$
    \begin{pmatrix}f_1\\\vdots\\ f_t\end{pmatrix}\big|U_{\ell^2}^2
   =(A_1^2+C_1)\begin{pmatrix}f_1\\\vdots\\ f_t\end{pmatrix}
   +A_1B_1\begin{pmatrix}g_1\\\vdots\\ g_t\end{pmatrix}
   +A_1C_1\begin{pmatrix}f_1\\\vdots\\
     f_t\end{pmatrix}\big|V_{\ell^2}.
  $$
  Iterating, we see that in general if
  $$
    \begin{pmatrix}f_1\\\vdots\\ f_t\end{pmatrix}\big|U_{\ell^2}^k
   =A_k\begin{pmatrix}f_1\\\vdots\\ f_t\end{pmatrix}
   +B_k\begin{pmatrix}g_1\\\vdots\\ g_t\end{pmatrix}
   +C_k\begin{pmatrix}f_1\\\vdots\\
     f_t\end{pmatrix}\big|V_{\ell^2},
  $$
  then the coefficients satisfy the recursive relation
  $$
    A_{k+1}=A_kA_1+C_k, \qquad B_{k+1}=A_kB_1, \qquad C_{k+1}=A_kC_1.
  $$
  (Note that $B_1$ and $C_1$ are scalar matrices. Thus, all
  coefficients are polynomials in $A$.) Finally, we note that the
  relation $A_{k+1}=A_kA_1+C_k=A_kA_1+C_1A_{k-1}$ can be written as
  $$
    \begin{pmatrix}A_{k+1}\\A_k\end{pmatrix}
   =\begin{pmatrix}A&C_1\\I_t&0\end{pmatrix}
    \begin{pmatrix}A_k\\A_{k-1}\end{pmatrix}.
  $$
  which yields
  $$
    \begin{pmatrix}A_{k+1}\\A_k\end{pmatrix}
   =\begin{pmatrix}A&C_1\\I_t&0\end{pmatrix}^k
    \begin{pmatrix}A\\I_t\end{pmatrix}
   =\begin{pmatrix}A&C_1\\I_t&0\end{pmatrix}^{k+1}
    \begin{pmatrix}I_t\\0\end{pmatrix}.
  $$
  This proves the theorem.
\end{proof}

\begin{proof}[Proof of Theorem \ref{theorem: congruences}]
  Let $m\ge 13$ be a prime. Let $r$ be the integer satisfying
  $0<r<24$ and $m\equiv-r\mod 24$ and set $s=(m-r-2)/2$. By Corollary
  \ref{corollary: Chua}, $F_{m,1}$ congruent to a modular form in
  $\S_{r,s}$, where $\S_{r,s}$ is defined by \eqref{equation: Srs}.
  Now let $\{f_1,\ldots,f_t\}$ be a basis for $\S_{r,s}$ and $A$ be
  given as in the statement of the theorem. Then by Corollary
  \ref{corollary: recursion}, we know that
  $$
    \begin{pmatrix}f_1\\\vdots\\f_t\end{pmatrix}\big|U_{\ell^2}^k
=A_k\begin{pmatrix}f_1\\\vdots\\f_t\end{pmatrix}
+B_k\begin{pmatrix}g_1\\\vdots\\g_t\end{pmatrix}
+C_k\begin{pmatrix}f_1\\\vdots\\f_t\end{pmatrix}\big|V_{\ell^2},
  $$
  where $g_j=f_j\otimes \left(\frac\cdot\ell\right)$, and $A_k$,
  $B_k$, and $C_k$ are $t\times t$ matrices satisfying
%  the recursive relations
%  \begin{equation} \label{equation: A1}
%    A_1=A, \quad
%    B_1=-\ell^{(m-5)/2}\left(\frac{(-1)^{(r-1)/2}12}\ell\right)I_t,
%    \quad C_1=-\ell^{m-4}I_t,
%  \end{equation}
%  and
%  \begin{equation} \label{equation: Ak}
%    A_k=A_{k-1}A_1+C_{k-1}, \qquad B_k=A_{k-1}B_1, \qquad
%    C_k=A_{k-1}C_1.
%  \end{equation}
%  Set $A_0=I_t$ and
%  $$
%    X=\begin{pmatrix} A&-\ell^{m-4}I_t\\I_t&0\end{pmatrix}.
%  $$
%  Then the recursive relation for $A_k$ can be written as
  \begin{equation} \label{equation: recursion in matrices}
    \begin{pmatrix} A_k\\A_{k-1}\end{pmatrix}
%   =X\begin{pmatrix} A_{k-1}\\ A_{k-2}\end{pmatrix}
%   =\ldots
%   =X^{k-1}\begin{pmatrix} A_1\\A_0\end{pmatrix}
   =X^k\begin{pmatrix}I_t\\0\end{pmatrix},
  \end{equation}
  \begin{equation} \label{equation: recursion in matrices 2}
    B_k=-\ell^{(m-5)/2}\left(\frac{(-1)^{(r-1)/2}12}\ell\right)A_{k-1},
    \qquad
    C_k=-\ell^{m-4}A_{k-1}
  \end{equation}
  with
  $$
    X=\begin{pmatrix} A&-\ell^{m-4}I_t\\I_t&0\end{pmatrix}
  $$
  for all $k\ge 1$. Now we have
  $$
    X^{-1}=\ell^{-(m-4)}\begin{pmatrix}0&\ell^{m-4}I_t\\
    -I_t&A\end{pmatrix}
%   =\ell^{-2(m-4)}\begin{pmatrix}-\ell^{m-4}I_t&\ell^{m-4}A\\
%   -A&A^2-\ell^{m-4}I_t\end{pmatrix}.
  $$
  Therefore, if the order of $X\mod m$ in $\mathrm{PGL}(2t,\F_m)$ is
  $K$, then we have
  $$
    \begin{pmatrix}A_{uK-1}\\A_{uK-2}\end{pmatrix}
   =X^{uK-1}\begin{pmatrix}I_t\\0\end{pmatrix}
    \equiv\begin{pmatrix} 0\\U\end{pmatrix} \mod m.
  $$
  for some $t\times t$ matrix $U$, that is, $A_{uK-1}\equiv 0\mod m$.
  The rest of proof follows Ono's argument.

  We have
  $$
    \begin{pmatrix}f_1\\\vdots\\f_t\end{pmatrix}\big|U_{\ell^2}^{uK-1}
    \equiv B_{uK-1}\begin{pmatrix}g_1\\\vdots\\g_t\end{pmatrix}
   +C_{uK-1}\begin{pmatrix}f_1\\\vdots\\f_t\end{pmatrix}\big|V_{\ell^2}
    \mod m
  $$
  and
  $$
    \begin{pmatrix}f_1\\\vdots\\f_t\end{pmatrix}\big|U_{\ell^2}^{uK-1}
    \big|U_\ell\equiv C_{uK-1}
    \begin{pmatrix}f_1\\\vdots\\f_t\end{pmatrix}\big|V_{\ell}\mod m.
  $$
  This implies that the $\ell^{2uK-1}n$th Fourier coefficients of $f_j$
  vanishes modulo $m$ for all $j$ and all $n$ not divisible by $\ell$.
  Since $F_{m,1}$ is a linear combination of $f_j$ modulo $m$, the
  same thing is true for the $\ell^{2uK-1}n$th Fourier coefficients of
  $F_{m,1}$. This translates to
  $$
    p\left(\frac{m\ell^{2uK-1}n+1}{24}\right)\equiv 0\mod m
  $$
  for all $n$ not divisible by $\ell$. This proves \eqref{equation:
    result} for the case $j=1$.

  For the case $j>1$, we note that the operators $U_\ell$ and $U_m$
  commutes. Thus,
  $$
    \begin{pmatrix}f_1\\\vdots\\f_t\end{pmatrix}\big|U_m^j\big|U_{\ell^2}^k
   =A_k\begin{pmatrix}f_1\\\vdots\\f_t\end{pmatrix}\big|U_m^j
   +B_k\begin{pmatrix}g_1\\\vdots\\g_t\end{pmatrix}\big|U_m^j
   +C_k\begin{pmatrix}f_1\\\vdots\\f_t\end{pmatrix}\big|V_{\ell^2}
    \big|U_m^j,
  $$
  where $A_k$, $B_k$, and $C_k$ satisfy the same relations
  \eqref{equation: recursion in matrices} and \eqref{equation:
    recursion in matrices 2}. Taking the fact
  \eqref{equation: F|Um} into account, we see that the same argument
  in the case $j=1$ gives us the general congruence.

  Finally, if the matrix $X$ has order $M$ in $\mathrm{GL}(2t,\F_m)$,
  then from the recursive relations \eqref{equation: recursion in
    matrices} and \eqref{equation: recursion in matrices 2},
  it is obvious that \eqref{equation: periodicity result} holds.
  This completes the proof.
\end{proof}

\begin{Remark} In general, the integer $K$ in Theorem \ref{theorem:
    congruences} may not be the smallest positive integer such that
  congruence \eqref{equation: Yang congruences} hold. To see this, for
  simplicity, we assume that $\S_{r,s}$ has dimension $t$ and its
  reduction modulo $m$ has a basis consisting of Hecke eigenforms
  $f_1,\ldots,f_t$ defined over $\F_m$. If the eigenvalues of
  $T_{\ell^2}$ for $f_i$ modulo $m$ are
  $a_\ell^{(1)},\ldots,a_\ell^{(t)}\in\F_m$. Let $k_i$ denote
  the order of
  $\left(\begin{smallmatrix}a_\ell^{(i)}&-\ell^{m-4}\\1&0\end{smallmatrix}
   \right)$ in $\mathrm{PGL}(2,\F_m)$. Let $k$ be the least
   common multiple of $k_i$. Then we can show that
  $$
    f_i\big|U_\ell^{2k-1}\equiv c_if_i\big|V_\ell\mod m
  $$
  for some $c_i\in\F_m$ and consequently congruence \eqref{equation:
    Yang congruences} holds. Of course, the least common multiple of
  $k_i$ may be smaller than the integer $K$ in Theorem \ref{theorem:
    congruences} in general.
\end{Remark}
\end{section}

\begin{section}{Examples} \label{section: examples}
\begin{Example} Let $m=13$. According to Corollary \ref{corollary:
  Chua}, we have
  $$
    F_{13,1}\equiv c\eta(24\tau)^{11} \mod 13
  $$
  for some $c\in\F_{13}$. (In fact, $c=11$. See \cite[page
  303]{Ono-Annals}.) The eigenvalues $a_\ell$ modulo $13$ of
  $T_{\ell^2}$ for the first few primes $\ell$ are
  $$ \extrarowheight3pt
  \begin{array}{c|ccccccccccccccccc} \hline\hline
  \ell & 5 & 7 & 11 & 17 & 19 & 23 & 29 & 31 & 37 & 41 & 43 & 47 & 53
       & 59 & 61 & 67 & 73 \\ \hline
a_\ell &10 & 8 &  5 & 1  &  8 & 8  &  4 &  4 &  5 &  9 & 12 &  6 & 10
       &  0 &  2 &  4 & 0 \\
\ell^9 & 5 & 8 & 8  & 12 &  5 & 12 &  1 &  5 &  8 &  5 & 12 &  8 &  1
       &  8 &  1 &  5 & 5 \\ \hline\hline
  \end{array}
  $$
  For $\ell=5$, the matrix
  $$
    X=\begin{pmatrix}a_\ell&-\ell^9\\1&0\end{pmatrix}
      \equiv\begin{pmatrix}10&8\\1&0\end{pmatrix} \mod 13
  $$
  has eigenvalues $5\pm\sqrt 7$ over $\F_{13}$. Now the order of
  $(5+\sqrt7)/(5-\sqrt 7)$ in $\F_{169}$ is $14$. This implies that
  $14$ is the order of $X$ in $\mathrm{PGL}(2,\F_{13})$ and we have
  $$
    p\left(\frac{13\cdot5^{28u-1}n+1}{24}\right)\equiv 0\mod 13
  $$
  for all positive integers $u$ and all positive integers $n$ not
  divisible by $5$. Likewise, we find congruence \eqref{equation: Yang
    congruences} holds with
  $$ \extrarowheight3pt
  \begin{array}{c|ccccccccccccccccc} \hline\hline
\ell & 5 & 7 & 11 & 17 & 19 & 23 & 29 & 31 & 37 & 41 & 43 & 47
     & 53 & 59 & 61 & 67 & 73 \\ \hline
k    &14 &14 & 14 & 7  & 14 & 3  &  6 & 12 & 14 & 12 &  7 & 12
     & 7  &  2 & 13 & 12 & 2  \\
\hline\hline
  \end{array}
  $$
\end{Example}

\begin{Example}
  Let $m=37$. By Corollary \ref{corollary: Chua}, we know that
  $F_{37,1}$ is congruent to a cusp form in $\S_{11,12}$ modulo $37$.
  In fact, according to \cite[Table 3.1]{Chua-Archiv},
  $$
    F_{37,1}\equiv\eta(24\tau)^{11}(E_4(24\tau)^3+17\Delta(24\tau))
    \mod 37.
  $$
  The two eigenforms of $\S_{11,12}$ are defined over a certain
  real quadratic number field, but the reduction of
  $\S_{11,12}\cap\Z[[q]]$ modulo $37$ has eigenforms defined over
  $\F_{37}$. They are
  $$
    f_1=\eta(24\tau)^{11}(E_4(24\tau)^3+24\Delta(24\tau)),
    \qquad f_2=\eta(24\tau)^{11}\Delta(24\tau).
  $$
  Let $a_\ell^{(i)}$ denote the eigenvalue of $T_{\ell^2}$ associated
  to $f_i$. We have the following data.
  $$ \extrarowheight3pt
  \begin{array}{c|ccccccccccccccc} \hline\hline
\ell & 5 & 7 & 11 & 13 & 17 & 19 & 23 & 29 & 31 & 41
     &43 &47 & 53 & 59 & 61 \\ \hline
a_\ell^{(1)} & 1 & 33 & 22 & 7 & 11 & 0 & 1 & 9 & 35 & 11
     & 28 & 14 & 30 & 24 & 12 \\
a_\ell^{(2)} &32 & 10 &  0 & 6 &  7 & 8 & 31&36 & 9  & 10
     & 1  & 35 &  9 &  3 & 16 \\
\ell^{33} &  8 & 26 & 36 &  8 & 23 & 8 & 6 & 31 & 31 & 11
     & 6 & 1 & 10 & 23 & 29 \\ \hline\hline
  \end{array}
  $$
  Let
  $$
    X_i=\begin{pmatrix}a_\ell^{(i)} & -\ell^{33} \\ 1 &
      0\end{pmatrix}.
  $$
  For $\ell=5$, we find the orders of $X_1$ and $X_2$ in
  $\mathrm{PGL}(2,\F_{37})$ are $38$ and $12$, respectively. The least
  common multiple of the orders is $228$. Thus, we have
  $$
    p\left(\frac{37\cdot5^{456u-1}n+1}{24}\right)\equiv 0\mod 37
  $$
  for all positive integers $u$ and all positive integers $n$ not
  divisible by $5$. Note that this is an example showing that the
  integer $K$ in the statement of Theorem \ref{theorem: congruences}
  is not optimal. (Here we have $K=456$.)

  For other small primes $\ell$, we find congruence
  $$
    p\left(\frac{37\ell^{2uk-1}n+1}{24}\right)\equiv 0\mod 37
  $$
  holds for all $n$ not divisible by $\ell$ with
  $$ \extrarowheight3pt
  \begin{array}{c|ccccccccccccccc} \hline\hline
\ell & 5 & 7 & 11 & 13 & 17 & 19 & 23 & 29 & 31 & 41
     &43 &47 & 53 & 59 & 61 \\ \hline
k    &228&57 & 18 &684 & 38 & 38 &684 &684 &228 &171
     &18 &333& 18 & 12 & 684 \\ \hline\hline
  \end{array}
  $$
\end{Example}
\end{section}

\begin{section}{Generalizations} \label{section: generalizations}
  There are several directions one may
  generalize Theorem \ref{theorem: congruences}. Here we only consider
  congruences of the partition function modulo prime powers. The case
  $m=5$ will be dealt with separately because in this case we have
  a very precise congruence result.
% Later on, we will also consider
%  congruences for powers of the partition function.

  In his proof of Ramanujan's conjecture for the cases $m=5,7$, Watson
  \cite[page 111]{Watson-Crelle} established a formula
  $$
    F_{5,j}=\begin{cases}
    \displaystyle
    \sum_{i\ge 1}c_{j,i}\frac{\eta(120\tau)^{6i-1}}{\eta(24\tau)^{6i}},
   &\text{if }j\text{ is odd}, \\
    \displaystyle
    \sum_{i\ge 1}c_{j,i}\frac{\eta(120\tau)^{6i}}{\eta(24\tau)^{6i+1}},
    &\text{if }j\text{ is even},
    \end{cases}
  $$
  where
  $$
    c_{j,i}\equiv\begin{cases}3^{j-1}5^j \mod 5^{j+1}, &\text{if }i=1, \\
    0 \mod 5^{j+1}, &\text{if }i\ge 2. \end{cases}
  $$
  From the identity, one deduces that
  \begin{equation} \label{equation: modulo power of 5}
    F_{5,j}\equiv 3^{j-1}5^j\begin{cases}
    \eta(24\tau)^{19} \mod 5^{j+1}, &\text{if }j\text{ is odd}, \\
    \eta(24\tau)^{23} \mod 5^{j+1}, &\text{if }j\text{ is even}.
    \end{cases}
  \end{equation}
  Then Lovejoy and Ono \cite{Lovejoy-Ono-Crelle} used this formula to
  study congruences of the partition function modulo higher powers of
  $5$. One distinct feature of \cite{Lovejoy-Ono-Crelle} is the
  following lemma.

  \begin{Lemma}[{Lovejoy and Ono \cite[Theorem 2.2]{Lovejoy-Ono-Crelle}}]
    \label{lemma: Lovejoy-Ono}
    Let $\ell\ge 5$ be a prime. Let $a$ and $b$ be the
    eigenvalues of $\eta(24\tau)^{19}$ and $\eta(24\tau)^{23}$ for the
    Hecke operator $T_{\ell^2}$, respectively. Then we have
    $$
      a,b\equiv\left(\frac{15}\ell\right)(1+\ell) \mod 5.
    $$
  \end{Lemma}

  With this lemma, Lovejoy and Ono obtained congruences of the form
  $$
    p\left(\frac{5^j\ell^kn+1}{24}\right)\equiv 0\mod 5^{j+1}
  $$
  for primes $\ell$ congruent to $3$ or $4$ modulo $5$. Here we
  shall deduce new congruences using our method.

  \begin{Theorem} \label{theorem: power of 5}
  Let $\ell\ge 7$ be a prime. Set
  $$
    K_\ell=\begin{cases}
    5, &\text{if }\ell\equiv 1\mod 5, \\
    4, &\text{if }\ell\equiv 2,3\mod 5, \\
    2, &\text{if }\ell\equiv 4\mod 5.
    \end{cases}
  $$
  Then we have
  $$
    p\left(\frac{5^j\ell^{2uK_\ell-1}n+1}{24}\right)\equiv 0\mod 5^{j+1}
  $$
  for all positive integers $j$ and $u$ and all integers $n$ not
  divisible by $\ell$.
  \end{Theorem}

  \begin{proof} In view of \eqref{equation: modulo power of 5}, We
    need to study when a Fourier coefficient of $\eta(24\tau)^{19}$ or
    $\eta(24\tau)^{23}$ vanishes modulo $5$.

    Let $f=\eta(24\tau)^{19}$. Let $\ell\ge 7$ be a prime and $a$ be
    the eigenvalue of $T_{\ell^2}$ associated to $f$. By Corollary
    \ref{corollary: recursion} we have
    \begin{equation} \label{equation: mod 5}
      f\big|U_{\ell^2}^k=a_kf+b_kf\otimes\left(\frac\cdot\ell\right)
      +c_kf\big|V_{\ell^2},
    \end{equation}
    where $a_1=a$, $b_1=-\ell^8\left(\frac{-12}\ell\right)$,
    $c_1=-\ell^{17}$, and $a_k=a_{k-1}a_1+c_{k-1}$, $b_k=a_{k-1}b_1$,
    $c_k=a_{k-1}c_1$. According to the proof of Theorem \ref{theorem:
      congruences}, if the order of
    \begin{equation} \label{equation: temp theorem 8}
      \begin{pmatrix}a&-\ell^{17}\\1&0\end{pmatrix} \mod5
    \end{equation}
    in $\mathrm{PGL}(\F_5)$ is $k$, then
    \begin{equation} \label{equation: temp 2 theorem 8}
      f\big|U_{\ell^{2uk-1}}\equiv f\big|V_\ell \mod 5
    \end{equation}
    for all positive integers $u$. Now by Lemma \ref{lemma:
      Lovejoy-Ono} the characteristic polynomial of \eqref{equation:
      temp theorem 8} has a factorization
    $$
      \left(x-\left(\frac{15}\ell\right)\right)
      \left(x-\left(\frac{15}\ell\right)\ell\right)
    $$
    modulo $5$. From this we see that the order of \eqref{equation:
      temp theorem 8} in $\mathrm{PGL}(\F_5)$ is
    $$
      K_\ell=\begin{cases}
      5, &\text{if }\ell\equiv 1\mod 5, \\
      4, &\text{if }\ell\equiv 2,3\mod 5, \\
      2, &\text{if }\ell\equiv 4\mod 5. \end{cases}
    $$
    Thus, \eqref{equation: temp 2 theorem 8} holds with $k=K_\ell$.
    This yields the congruence
    $$
      p\left(\frac{5^j\ell^{2uK_\ell-1}n+1}{24}\right)\equiv 0\mod 5^{j+1}
    $$
    for odd $j$, positive integer $u$, and all positive integers $n$
    not divisible by $\ell$.

    The proof of the case $j$ even is similar to the above
    and is omitted.
  \end{proof}

  \begin{Remark} In \cite{Watson-Crelle}, Watson also had an identity
    for $F_{7,j}$, with which one can study congruences modulo higher
    powers of $7$. However, because there does not seem to exist an
    analog of Lemma \ref{lemma: Lovejoy-Ono} in this case, we do not
    have a result as precise as Theorem \ref{theorem: power of 5}
  \end{Remark}

  The next congruence result is an analog of Theorem 2 of
  \cite{Weaver-RamanujanJ}, which in turn is originated from the
  argument outlined in \cite[page 301]{Ono-Annals}.

  \begin{Theorem} Let $\ell\ge 7$ be a prime. Assume that one of the
    following situations occurs.
  \begin{enumerate}
  \item $\ell\equiv 1\mod 5$, $\left(\frac{-n}\ell\right)=-1$ with
    $k_\ell=2$ and $m_\ell=5$, or
  \item $\ell\equiv 2\mod 5$, $\left(\frac{-n}\ell\right)=(-1)^{i-1}$
    with $k_\ell=2$ and $m_\ell=4$, or
  \item $\ell\equiv 3\mod 5$, $\left(\frac{-n}\ell\right)=(-1)^{i-1}$
    with $k_\ell=1$ and $m_\ell=4$.
  \end{enumerate}
  Then
  $$
    p\left(\frac{5^i\ell^{2(um_\ell+k_\ell)}n+1}{24}\right)\equiv 0\mod 5^{i+1}
  $$
  for all non-negative integers $u$.
  \end{Theorem}

  \begin{proof} Assume first that $i$ is odd. Again, in view of
    \eqref{equation: modulo power of 5}, we need to study when the
    Fourier coefficients of $f(\tau)=\eta(24\tau)^{19}$ vanish modulo
    $5$.

  Let $\ell\ge 7$ be a prime and $a$ be the eigenvalue of $T_{\ell^2}$
  associated to $\ell$. By \eqref{equation: mod 5}, we have
  \begin{equation} \label{equation: mod 5 2}
    f\big|U_{\ell^2}^k=a_kf+b_kf\otimes\left(\frac\cdot\ell\right)
   +c_kf\big|V_{\ell^2},
  \end{equation}
  where $a_k$, $b_k$, $c_k$ satisfy
  $$
    \begin{pmatrix}a_k\\a_{k-1}\end{pmatrix}
   =\begin{pmatrix}a&-\ell^{17}\\1&0\end{pmatrix}^k
    \begin{pmatrix}1\\0\end{pmatrix}, \quad
    b_k\equiv-\left(\frac{-12}\ell\right)a_{k-1},  \quad
    c_k\equiv-\ell a_{k-1} \mod 5.
  $$
  From Lemma \ref{lemma: Lovejoy-Ono}, we know that for $\ell\equiv
  1\mod5$, we have $a_1\equiv 2\epsilon$ and thus the
  values of $a_k$ modulo $5$ are
  $$ \extrarowheight3pt
  \begin{array}{ccccccccccccc} \hline\hline
a_1 & a_2 & a_3 & a_4 & a_5 & a_6 & a_7 & a_8 & a_9 & a_{10} & a_{11}
& a_{12} & \ldots \\ \hline
2\epsilon & 3 & 4\epsilon & 0 & \epsilon & 2 & 3\epsilon & 4 & 0 & 1
& 2\epsilon & 3 & \ldots \\ \hline\hline
  \end{array}
  $$
  where $\epsilon=\left(\frac{15}\ell\right)$. Now assume that
  $f(\tau)=\sum c(n)q^n$. Comparing the $n$th Fourier coefficients of
  the two sides of \eqref{equation: mod 5 2} for integers $n$
  relatively prime to $\ell$, we obtain
  $$
    c(\ell^{2k}n)=\left(a_k+b_k\left(\frac n\ell\right)\right)c(n)
    \equiv\left(a_k-a_{k-1}\left(\frac{-12n}\ell\right)\right)c(n)\mod 5.
  $$
  When $k=5u+2$ for a non-negative integer $u$, we have
  \begin{equation} \label{equation: mod 5 3}
  \begin{split}
    c(\ell^{2(5u+2)}n)&\equiv 3\left(\frac{15}\ell\right)^u\left(
    1+\left(\frac{15}\ell\right)\left(\frac{-12n}\ell\right)\right)c(n) \\
  &=3\left(\frac{15}\ell\right)^u
    \left(1+\left(\frac{-n}\ell\right)\right)c(n) \mod 5.
  \end{split}
  \end{equation}
  Thus, if $\left(\frac{-n}\ell\right)=-1$, then
  $c(\ell^{2(5u+2)}n)\equiv 0\mod 5$. This translates to the
  congruence
  $$
    p\left(\frac{5^i\ell^{2(5u+2)}n+1}{24}\right)\equiv 0\mod 5^{i+1}.
  $$
  This proves the first case of the theorem. The proof of the other
  cases is similar.
  \end{proof}

  \begin{Example}
  \begin{enumerate}
  \item Let $\ell=11$, $i=1$, and $n=67$. Then the first situation
    occurs. We find
    $$
      p\left(\frac{5\cdot11^4\cdot67+1}{24}\right)
     =p(204364)=28469\ldots\ldots\ldots24450,
    $$
    which is a multiple of $25$.
  \item Let $\ell=11$, $i=1$, and $n=19$. The condition in the theorem
    is not fulfilled, but \eqref{equation: mod 5 3} implies that
  $$
    p\left(\frac{5\cdot11^4\cdot19+1}{24}\right)
    \equiv p\left(\frac{5\cdot19+1}{24}\right) \mod 25.
  $$
  Indeed, we have $p(4)=5$,
  $$
    p(57954)=37834\ldots\ldots\ldots45055,
  $$
  and they are congruent to each other modulo $25$.
  \item Let $\ell=7$, $i=2$, and $n=23$. Then the second
    situation occurs. We have
    $$
      p\left(\frac{5^2\cdot7^4\cdot23+1}{24}\right)
     =p(57524)=38402\ldots\ldots\ldots43875,
    $$
    which is indeed a multiple of $5^3$.
  \end{enumerate}
  \end{Example}

  \begin{Theorem} \label{theorem: prime powers} Let $m\ge 13$ be a
    prime and $\ell$ be a prime different from $2,3,m$. For each
    positive integer $i$, there exists a positive integer $K$ such
    that for all $j\ge i$, all $u\ge 1$ and all positive integers $n$
    not divisible by $\ell$, the congruence
    $$
      p\left(\frac{m^j\ell^{2uK-1}n+1}{24}\right)\equiv 0\mod m^i
    $$
    holds. There is also another positive integer $M$ such that
    $$
      p\left(\frac{m^j\ell^rn+1}{24}\right)\equiv
      p\left(\frac{m^j\ell^{M+r}n+1}{24}\right)\mod m^i
    $$
    holds for all $n$.
  \end{Theorem}

  \begin{proof} Let $\beta_{m,j}$ be the integer satisfying
    $1\le\beta_{m,i}\le m^i-1$ and $24\beta_{m,i}\equiv 1\mod m^i$.
    Define
    $$
      k_{m,i}=\begin{cases}
      (m^{i-1}+1)(m-1)/2-12\lfloor m/24\rfloor-12,
      &\text{if }i\text{ is odd}, \\
      m^{i-1}(m-1)-12, &\text{if }i\text{ is even}.
      \end{cases}
    $$
    By Theorem 3 of \cite{Ahlgren-Boylan-Inventiones}, for all $i\ge
    1$, there is a modular form $f\in M_{k_{m,i}}(\SL(2,\Z))$ such that
    $$
      F_{m,i}\equiv \eta(24\tau)^{(24\beta_{m,i}-1)/m^i}
      f(24\tau) \mod m^i.
    $$
    The rest of proof is parallel to that of Theorem \ref{theorem:
      congruences}.
  \end{proof}

  \begin{Example} Consider the case $m=13$ and $i=2$ of Theorem
    \ref{theorem: prime powers} and assume that $\ell$ is a prime
    different from $2,3,13$. By \cite[Theorem
    3]{Ahlgren-Boylan-Inventiones}, $F_{13,2}$ is congruent to a
    modular form in the space $\S_{23,144}$ of dimension $13$. Choose
    a $\Z$-basis
    $$
      f_i=\eta(24\tau)^{23}E_4(24\tau)^{3(13-i)}\Delta(24\tau)^{i-1}, \quad
      i=1,\ldots,13,
    $$
    for $\Z[[q]]\cap\S_{23,144}$ and let $A$ be the
    matrix of $T_{\ell^2}$ with respect to this basis. If the
    order of the matrix
    $$
      \begin{pmatrix}A&-\ell^{309}I_{13}\\I_{13}&0\end{pmatrix}
      \mod 169
    $$
    in $\mathrm{PGL}(26,\Z/169)$ is $K$, then we have
    $$
      p\left(\frac{169\ell^{2K-1}n+1}{24}\right)\equiv 0\mod 169
    $$
    for all integers $n$ not divisible by $\ell$.
    For instance, for $\ell=5$, we find
  $$ \tiny A=\left(
  \begin{array}{ccccccccccccc}
20& 101& 52& 52& 166& 148& 46& 135& 96& 51& 73& 49& 128 \\
166& 164& 159& 66& 123& 50& 144& 85& 29& 116& 22& 93& 10 \\
158& 152& 90& 65& 20& 167& 27& 96& 109& 154& 127& 164& 76 \\
120& 154& 132& 110& 22& 113& 115& 51& 25& 104& 108& 82& 33 \\
43& 148& 131& 45& 81& 2& 164& 145& 117& 157& 4& 108& 61 \\
134& 23& 151& 120& 151& 44& 30& 1& 76& 32& 60& 132& 165 \\
121& 40& 83& 4& 56& 88& 3& 134& 100& 85& 88& 18& 3 \\
23& 20& 20& 31& 66& 24& 41& 126& 47& 137& 33& 112& 49 \\
143& 18& 44& 26& 89& 109& 118& 148& 35& 16& 35& 122& 150 \\
144& 51& 47& 143& 109& 164& 52& 38& 92& 50& 98& 60& 104 \\
70& 165& 89& 80& 28& 75& 19& 110& 101& 41& 155& 78& 67 \\
123& 147& 54& 4& 60& 133& 49& 151& 30& 32& 157& 108& 82 \\
95& 139& 50& 70& 124& 168& 87& 63& 13& 104& 58& 107& 113
  \end{array} \right)
  $$
  modulo $169$, and the order $K$ is $28392$, which yields
  $$
    p\left(\frac{13^2\cdot5^{56783}n+1}{24}\right)\equiv 0\mod 13^2
  $$
  for all $n$ not divisible by $5$.
%    It turns out that in this case, we can do a little better.
%    Over the $13$-adic field $\Q_{13}$, the Hecke eigenforms of
%    $\S_{23,144}$ are all defined over $\Q_{13}$. For example, if we
%    let $f_r=\eta(24\tau)^{23}E_4(24\tau)^{36-3r}\Delta(24\tau)^r$,
%    $r=0,\ldots,12$, then one of the Hecke eigenform is
%    \begin{equation*}
%    \begin{split}
%      g_1&=f_0+131698f_1+167534f_2+240473f_3+212749f_4+175490f_5\\
%       &\qquad\qquad+135582f_6+36904f_7+235050f_8+53606f_9+162776f_{10}\\
%       &\qquad\qquad+30625f_{11}+360741f_{12}+O(13^5)
%    \end{split}
%    \end{equation*}
%    and its eigenvalue for $T_{25}$ is $a_5=130843+O(13^5)$. The matrix
%    \begin{equation} \label{equation: example prime power}
%      \begin{pmatrix}a_5&-5^{309}\\1&0\end{pmatrix} \mod 169
%    \end{equation}
%    has order $156$ in $\mathrm{PGL}(2,\Z/169)$.
%    For the other $12$ Hecke eigenforms $g_r$, $r=2,\ldots,13$, we
%    find the eigenvalues modulo $169$ for $T_{25}$ are $145$, $132$,
%    $11$, $2$, $50$, $10$, $54$, $75$, $166$, $63$, $58$, and $127$
%    with the orders of the corresponding matrices \eqref{equation:
%      example prime power} being one of $156$, $182$, and $52$. Thus,
%    $$
%      g_r\big|U_5^{2184u-1}\equiv g_r\big|V_5 \mod 169
%    $$
%    for all $r$ and we have
%    $$
%      p\left(\frac{13^2\cdot5^{2183}n+1}{24}\right)\equiv 0\mod 169
%    $$
%    for all positive integers $n$ not divisible by $5$.
  \end{Example}

%  In the remainder of the section, we will consider congruences for
%  \emph{powers} of the partition function, which, for $r\in\Z$, are
%  defined in terms of the generating function
%  $$
%    \sum_{n=0}^\infty p_r(n)q^n=\prod_{n=1}^\infty\frac1{(1-q^n)^r}.
%  $$
%  (In particular, $p_1(n)$ is the usual partition function.) The
%  congruence properties of these functions have also been studied by
%  many authors. See
%  \cite{Atkin-CJM,Boylan-AA,Kiming-Olsson-Archiv,Newman-CJM,Stanger-RamanujanJ}
%  and the references contained therein.
\end{section}

\bibliographystyle{plain}

\end{document}